\documentclass[a4paper,10pt,twoside]{article}
\usepackage[french, english]{babel}
\usepackage[dvips]{graphics}
\usepackage{epsf}
\usepackage{graphicx}
\usepackage{epsfig}
\usepackage{amsmath}
\usepackage[psamsfonts]{amssymb}
\usepackage[psamsfonts]{eucal}
\usepackage{amsthm}
\usepackage[utf8]{inputenc}
\usepackage{textcomp}
\usepackage{mathptmx}
\usepackage[scaled]{helvet}
\usepackage{perso2}
\usepackage{xcolor}

\swapnumbers
\theoremstyle{plain}
\theoremstyle{definition}
\newcommand{\field}[1]{\mathbb{#1}}
\newcommand{\CC}{\field{C}}

\newcommand{\RR}{\field{R}}

\def\div{\mathop{\rm div}\nolimits}
\def\rot{\mathop{\rm rot}\nolimits}

\def\d{\mathrm{d}}
\def\mathi{\mathrm{i}}

\def\SO{\mathrm{SO}}
\def\SU{\mathrm{SU}}

\title{Presentation of Jean-Marie Souriau's book\\
        \lq\lq Structure des systèmes dynamiques\rq\rq}
\author{G\'ery de Saxc\'e\\
        Emeritus Professor\\
        Univ. Lille, CNRS, Centrale Lille, UMR 9013 – LaMcube – \\
        Laboratoire de m\'ecanique multiphysique multi\'echelle, F-59000, Lille, France\\\\
        Charles-Michel Marle\\
        Honorary Professor\\ retired from the former Université Pierre et Marie Curie
        \\ today Sorbonne Université, Paris, France
  }

\begin{document}

\maketitle

\begin{abstract}
Jean-Marie Souriau's book \lq\lq Structure des systèmes dynamiques\rq\rq, 
published in 1970, republished recently by Gabay, translated in English and 
published under the title  \lq\lq Structure of Dynamical Systems, a Symplectic 
View of Physics\rq\rq, is a work with an exceptional wealth which, 
fifty years after its publication, is still topical. In this paper, we give a 
rather detailled description of its content and we intend to highlight the 
ideas that to us, are the most creative and promising.
\end{abstract}
\vskip 1cm

Jean-Marie Souriau's book \lq\lq Structure des systèmes dynamiques\rq\rq\ 
\cite{SSD}  {was published in 1970 in a book collection for students  in the 
first year of master's degree in mathematics. It is directed in fact to mathematicians, 
beginners or experienced, wishing to know the applications of mathematics to 
physical sciences, and to physicists concerned with knowing certain mathematical 
tools useful for their researches.} The author was very aware of this since, 
in his Introduction,  he gives reading recommendations adapted to both reader 
categories. It is a work with an exceptional wealth which, fifty years after 
its publication, is still topical. In the first part, we shall be describing 
its content. In the second one, we shall discuss the most original aspects. 
As conclusion, we shall indicate why this book seem to us still today a valuable 
source for the students and researchers in mathematics, mechanics and physics. 

\section{Book content\label{Contenu}}

\subsection{A quick survey.}

This book  {comprises 
a rather large} introduction (20 pages) and five chapters. The first two chapters, 
entitled \emph{Differential geometry} and \emph{Symplectic geometry}, are purely mathematical. 
The third one, entitled \emph{mechanics}, begins with a classical presentation 
of the mechanics of material point systems. Very soon, the author introduces 
the concept of \emph{manifold of motions} of a  {system. Then he applies
the methods of symplectic geometry presented in the previous chapter both to
classical and relativistic mechanical systems}. An important paragraph deals with 
isolated systems admitting a symmetry group, acting transitively on the manifold 
of motions, that the author considers as models of \emph{elementary particles}. 
Next the author studies the dynamics of systems of elementary  {particles,} 
taking into account their interactions. The fourth chapter, entitled \emph{statistical mechanics}, 
{contains two sections.} 
The first one, essentially mathematical, presents  
{measure and integration theories,}
together with notions of probability theory. The author defines the \emph{statistical states} 
of a dynamical system, the \emph{entropy} of a statistical state and proposes a 
generalisation of the notion of \emph{Gibbs state}. In the second paragraph, 
he uses these notions to treat certain systems seen in physics: classical and 
relativistics ideal gas, systems of null mass particles, to name a few. 
He clarifies the interpretation of the parameters which the Gibbs state depends 
on in terms of thermodynamical quantities, and proposes an interesting 
generalisation of the concept of thermodynamical equilibrium. The fifth and 
last chapter, entitled \emph{Geometric quantization},  {presents a construction}  
which allows to associate to a symplectic manifold satisfying certain conditions 
another manifold, of  {larger dimension,  called by the author a} 
\emph{quantum manifold}. This construction is due to the author. A slightly 
 {different but equivalent construction} was proposed independently 
by the American  mathematician B.~Kostant \cite{Kostant70}. In the second and 
last  {
section} of the fifth chapter, the author applies this construction to the 
quantization of physical systems. 
\par\smallskip

The book is illustrated  {
by} many figures which mostly are very meaningful schematic representations of 
the geometric constructions used by the author. The references to other books or papers, 
fairly limited in number, are not gathered in a bibliographic list, but 
indicated in footnotes. An index, very detailled and easy to use, and the list 
of main notations completes the book.

\subsection{Detailed presentation}

\subsubsection*{The introduction.}
The author evocates the book by  Joseph-Louis
Lagrange (1736--1813) \emph{Mé\-ca\-ni\-que analytique} \cite{Lagrange5}, 
written at the end of the 18th century. This famous work is at the origin 
of the \emph{m\'ecanique analytique classique} which was, during the 19th century 
and the first half of the 20th century, an essential part of the scientific learning 
in the French universities and high schools. According to the author, it is an 
unfinished work in which some chapters are only sketched. 
For him, the form used to present this theory and the concepts it uses 
(instead of concepts, the author writes the \emph{categories} in the epistemological, 
Aristotelian or Kantian sense) were fixed by Lagrange's successors such as 
Sim\'eon Denis Poisson  (1781--1840), 
William Rowan Hamilton (1805--1865) and Carl Gustav Jakob Jacobi (1804--1851).
The author considers that the form thus given to  {
analytical mechanics}, although it gives to this theory a formal mathematical 
perfection, has lost an important part of Lagrange's thought. The great discoveries 
of the first quarter of the 20th century (special and general relativity, 
quantum theory) taught us  {that} the words \emph{time}, \emph{space} and 
\emph{matter} {not necessarily have}  the obvious meaning ascribed thereto. 
{For the author, classical analytical mechanics, which remains an essential 
ingredient of current physical theories, is not outdated, although certain concepts it uses 
are so because they do not have the required covariance}, in other words because 
they are in contradiction with Galilean relativity. 
He wishes to show in his book that a better consideration of Lagrange's thought 
allows a  formulation of this theory compatible with the most recent discoveries 
in physical sciences.
 \par\smallskip

 {Using the very simple example of the motion of a material point, 
the author explains which concepts he is going to use:} evolution space, 
space of motions of a system, Lagrange's form. The concept of  {space 
of motions of a dynamical system seems to us the most important~:} it is the set of  
{all possible motions of the considered system}. The author discusses 
 {in depth its usefulness, often underestimated by scholars mostly interested
in the study of one particular motion of the considered system.} 
He presents also all the mathematical tools that he will use: differential forms, 
Lie groups, symplectic forms, etc. 
Hence he reviews fairly  {accurately} every chapter of his book and 
completes his introduction with advises to the readers. 

\subsubsection*{The first chapter, Differential Geometry} 

This chapter presents in less than 70 pages numerous
delicate notions: differential manifolds, tangent and cotangent fiber bundles, 
submanifolds, covering spaces, vector fields and differential equations, 
Lie bracket of two vector fields, exterior derivative, foliations, Lie groups, 
calculus of variation. The author presents the concept of differential manifold 
in a rather original way  {which does not use the previous 
presentation of topological manifolds},
undoubtedly to make this notion easily accessible to beginning students. 
 {Differential manifolds are not assumed to be Hausdorff.  
The \emph{space of motions} of certain mechanical systems encountered in Chapter III
are indeed non-Hausdorff manifolds. The author's language sometimes slightly differs from
that generally used~: for example, he calls an \emph{embedding} what most of 
geometers call an \emph{injective immersion}. However, readers can easily avoid 
misunderstandings, since the author scrupulously defines all the terms he uses}.
 
The paragraph devoted to Lie groups contains a detailed presentation of the 
actions of a Lie group on a differential manifold and of its adjoint representation. 
The author will define the coadjoint representation in the next chapter, 
 {with the study of the moment map of the action of a dynamical group 
on a presymplectic or symplectic manifold}. The main classical Lie groups 
(linear, orthogonal, unitary, symplectic) are described 
in a very original manner. The section about 
calculus of variations presents, besides the classical Euler-Lagrange equations, 
the extremality criterion using Cartan's form, often called 
\emph{Euler-Cartan theorem}, that establishes a link between
calculus of variation and symplectic geometry, together with Noether's theorem.

\subsubsection*{Chapter II, Symplectic Geometry.} 
It is also essentially mathematical. Clearly shorter than the previous one (46 pages), 
it begins with the study of a  {finite dimensional vector space $E$ 
equipped with a bilinear skew-symmetric form $\sigma$. The author defines the 
concepts of \emph{orthogonality} with respect to $\sigma$, of \emph{isotropic} 
\emph{co-isotropic} and \emph{self-othogonal} vector subspaces, also called by other authors
\emph{Lagrangian vector subspaces} when $\sigma$ is nondegenerate}. 
The author proves  {that} the rank of  $\sigma$ is always even and
 {that given a coisotropic vector subspace}, one can always build 
a basis of the kernel of $\sigma$, then complete it to obtain a basis of $E$, 
called \emph{canonical}, in which the expression of $\sigma$ is very simple. 
 {When $\sigma$ is non degenerate, it is called a \emph{symplectic form}, 
the dimension of $E$ is even and the couple $(E,\sigma)$ is called a \emph{symplectic vector space}}. 
Hence the author defines the \emph{symplectic group} of $(E,\sigma)$ and studies 
its natural action on $E$, its features, together with  {thoses} 
of its Lie algebra. He shows in particular that $E$ has a complex vector space 
structure adapted, in a precise mathematical sense, to the symplectic 
form $\sigma$, and that it can be endowed with a Hermitian form of which  
$\sigma$ is the imaginary part.
\par\smallskip

The author defines next symplectic and presymplectic manifolds and studies 
their properties. He shows that under certain conditions, the quotient of a 
presymplectic manifold by the  {kernel  
of its presymplectic form} is a symplectic manifold, a result that he will use in the 
following chapter to define the \emph{space of motions} of a dynamical system. 
He shows that to each differentiable function defined on  {a} symplectic manifold, 
there is an associated vector field that he calls  {its}
\emph{symplectic gradient}.  {The flow of this vector field leaves unchanged 
the symplectic form}. He proves that the set of the differentiable functions 
defined on a symplectic manifold is endowed with a binary operation, the 
\emph{Poisson bracket}, that makes it a Lie algebra of finite dimension. He 
defines and studies some remarkable submanifolds of a symplectic manifold 
(isotropic, co-isotropic and self-orthogonal submanifolds), studies their 
properties and gives several examples thereof. He proves 
\emph{Darboux theorem}, whereby every point of a symplectic manifold is an element 
of the domain of a chart, called \emph{canonical}, in which the symplectic form 
is expressed in a simple manner. Next he defines and studies the 
\emph{symplectomorphisms}, also called \emph{canonical transformations}, 
and their generalisations, the  \emph{canonical similarities}, together with 
the \emph{infinitesimal canonical transformations} 
(called by other authors \emph{locally Hamiltonian vector fields}).
\par\smallskip

The author calls \emph{dynamical group} 
of a symplectic or presymplectic manifold a Lie group 
acting on it by canonical transformations. He calls \emph{moment} 
 {\emph{map}} of the dynamical 
group $G$ of a presymplectic or symplectic  manifold $M$ a smooth 
map $\Psi$ from $M$ into the vector space ${\mathcal G}^*$, dual of the Lie algebra 
$\mathcal G$ of $G$, such that for every $Z\in{\mathcal G}$, the infinitesimal 
generator of the action on $M$ of the one-parameter subgroup generated by  $Z$ 
is the symplectic gradient of the function linking, at every $x\in M$, the real 
$\bigl\langle \Psi(x),Z\bigr\rangle$. In the terminology used by most of 
geometers, this infinitesimal generator is  {called} 
the \emph{Hamiltonian vector field} 
whose Hamiltonian is the function  $x\mapsto\bigl\langle \Psi(x),Z\bigr\rangle$.

The author gives several examples of dynamical groups, indicates sufficient 
conditions for the existence of a moment map and studies its properties. 
This leads him to propose a generalisation of Noether's theorem encountered in 
the part of the previous chapter concerning the calculus of variations. 
After a quick presentation of the cohomology of Lie groups and algebras, 
the author shows that to any moment  {map} of the action of a dynamical group $G$ 
on a  {connected} presymplectic or 
symplectic manifold, there  {always exist} 
an associated cocycle  $\theta$ of $G$ valued in the dual ${\mathcal G}^*$ of its Lie algebra. 
The differential of $\theta$ at the neutral element, which is the cocycle
of the Lie algebra ${\mathcal G}$ associated to $\theta$, is a skew-symmetric
bilinear form on ${\mathcal G}$. 
The author calls \emph{symplectic cocycle} a 
cocycle of $G$ valued in ${\mathcal G}^*$ satisfying this property. He proves 
that the addition to the moment  {map} of a constant element of ${\mathcal G}^*$ 
modifies the cocycle $\theta$ by addition of a coboundary, therefore
does not  change its cohomology class, which depends only on the $G$-action, not 
on the choice of the moment  {map}. Moreover the cocycle $\theta$ allows to define an 
affine action of $G$ on the dual ${\mathcal G}^*$ of its Lie algebra  {whose} 
linear part is the coadjoint representation. The author then proves that 
the moment  {map} is equivariant with respect to 
 {this affine action of $G$ on ${\mathcal G}^*$ and its action on $M$
as a dynamical group of this manifold}. 
The orbits of this affine action are submanifolds embedded 
in  ${\mathcal G}^*$ (in author's meaning, most of the geometers would say 
rather \emph{immersed}) and are endowed with a symplectic form whose
expression involves the cocycle of $\mathcal G$, that is the differential of $\theta$ at the neutral element. 
Nowadays this important result is expressed by saying that ${\mathcal G}^*$ possesses 
a \emph{Poisson structure} whose
\emph{symplectic leaves} are the orbits of the affine action for which the moment  
 {map} is equivariant, that this
Poisson structure  {remains unchanged under this affine action}
and that  the moment map is a \emph{Poisson map}. When 
this Poisson structure on ${\mathcal G}^*$ 
was discovered by the Norwegian mathematician Sophus Lie (1842--1899) 
 {in the special case in which the cocycle $\theta$ vanishes}. 
In the general case, it was rediscovered independantly by  Alexander Kirillov, 
Bertram Kostant and the author, Jean-Marie Souriau. 
In his book \emph{Structure des systèmes dynamiques}, 
he does not use the concept of Poisson structure, rather he uses the fact that 
the orbits of the affine action of $G$ on ${\mathcal G}^*$ are symplectic manifolds, 
which he calls \emph{symplectic manifolds defined by a Lie group}. 
He shows that when the dynamical group $G$ of a symplectic manifold 
$(M,\sigma)$ is connected, and when  
its action is transitive and possesses a moment  {map} $\Psi$,
this map is a local symplectomorphism of $M$ onto 
an orbit of the affine action of  $G$ on ${\mathcal G}^*$ for which $\Psi$ is equivariant. 
The map $\Psi$ is a symplectomorphism if and only if, moreover, 
the isotropy group of a point of $M$ 
is connected. He presents some examples of this important result, often called nowadays 
\emph{Kostant-Souriau theorem}.

\subsubsection*{Chapter III, Mechanics.} 
Relatively long  (105 pages), 
this chapter begins with the study of a mechanical system composed of material 
points in a fixed Galilean frame. The author writes \emph{Newton's equation} 
expressing the equality of the force acting on each of those points and the 
product of the mass by its acceleration. He
treats briefly the case of a unique 
material point placed in a Coulomb field (Kepler's problem), next
the \emph{$N$-body problem}  {of} 
celestial mechanics. He introduces then the 
notion of \emph{constraint}, presents the \emph{principle of virtual work} and 
the conditions in which a constraint is called \emph{ideal}. He uses, to study 
the motion of a rigid body, the \emph{group of} 
 {\emph{Euclidean displacements}} 
and establishes the equations of motion. 
\par\smallskip

Returning to the case of a system of $N$ material points without constraints, 
subjected to forces expressed by differentiable functions of  {the}
time  {and of} the positions and  
velocities of these points, the author shows  {that} the equations 
of motion are expressed as 
the differential equation  associated to a 
vector field depending on time, defined on the \emph{evolution space} of the system 
( {a} set, isomorphic to $\RR^{6N+1}$,  made of
multiplets composed of the time and of the positions and velocities of the $N$ material points). These equations 
determine a  {\emph{foliation in curves}}
of the evolution space. Each of these curves 
is the mathematical expression of a possible motion of the system. 
It is why the author calls \emph{space of motions} of the system the set 
 {made} of these curves. He proves that the space of motions 
is a differential manifold (not always Hausdorff) of dimension $6N$ 
and there exists a natural projection, smooth
and anywhere of rank $6N$, of the evolution space on the space of motions. 
 {Then} he proves the existence, on the evolution space of the system, 
of a remarkable $2$-form $\sigma$, which he calls \emph{Lagrange's form}, 
because it was used  {in} 
1808 by Lagrange in his works in
celestial mechanics. As pointed out by the author, this form was used in mechanics
around 1950 by the French mathematician François Gallissot \cite{Gallissot}. 
For a unique material point of mass $m$ whose
position and velocity vectors are, respectively, ${\bf r}$ and ${\bf v}$, 
on which acts a force  ${\bf F}$, the components of these three vectors in a 
fixed orthonormal frame of the space being, respectively, 
$(r_1,r_2,r_3)$, $(v_1,v_2,v_3)$ and $(F_1,F_2,F_3)$,
this form reads
 $$\sigma=\sum_{i=1}^3(m\d v_i-F_i\d t)\wedge(\d r_i-v_i\d t)\,,\ 
    \hbox{that can be written}\ 
     \sigma= (m\d {\bf v}-{\bf F}\d t)\dot{\wedge}(\d{\bf r}-{\bf v}\d t)\,,
 $$
 where the symbol $\dot{\wedge}$ denotes the operator combining the dot product and the exterior product. 
 
The Lagrange form of a system of $N$ material points is the sum of the Lagrange 
forms of all points of the system. It determines the vector field of which the 
 {system's} motions are the integral curves, since its kernel is the sub-bundle 
that determines the foliation  {in} 
curves of the evolution space. The author remarks that it is always possible to 
choose $N$ differentiable vector fields  ${\bf B}_j$, defined on the evolution 
space, and to define $N$ other vector fields ${\bf E}_j$ so that  {the force} 
${\bf F}_j$ acting on the $j$-th material point is
${\bf F}_j={\bf E}_j- {\bf B}_j\times {\bf v}_j$, where ${\bf v}_j$ 
is the velocity of this material point. 
\par\smallskip

When changing the reference frame, the parameterization of the evolution space 
and the expression of the Lagrange form are modified, particularly because inertial 
forces must be included in the forces acting on the material points  
(centrifugal force and Coriolis' force).  The author states the
\emph{principe of Galilean relativity} which claims the existence of preferential 
reference frames, called \emph{inertial reference frames}, such that the 
expression of the Lagrange form of an \emph{isolated} system is the same in 
all inertial reference frames. He shows that the relative motion of an inertial 
reference frame with respect to another one is a translation motion at  {a}
constant velocity. The set of all changes of inertial reference frames is
a Lie group, called  {the} \emph{Galilei group}. 
The author gives a matrix expression of this group, which is 
of dimension $10$. It 
acts on the evolution space of an unconstrained system of material 
points by an action which preserves 
the Lagrange form. Using the principle of virtual works, the author extends 
the notions of evolution space,  
Lagrange form and space of motions, for 
a system of $N$ material points involving ideal 
constraints,  {which can be either holonomic or
non holonomic with a linear dependence on the velocities}. 
The results obtained for the unconstrained systems of material points 
remain  valid, provided that the constraints remain the same in all
inertial reference frames. 
\par\smallskip 
The author calls \emph{Maxwell's principle} the hypothesis that the exterior 
derivative of the Lagrange $2$-form of a general dynamical 
system (not necessarily made of material points) vanishes. For a system of material points, 
since the vector fields ${\bf B}_j$ which appear in the forces
${\bf F}_j={\bf E}_j-{\bf B}_j\times {\bf v}_j$ can be freely chosen,
the author formulates Maxwell's principle in the following form: the vector fields 
${\bf B}_j$ can be chosen in such a way that the Lagrange form is closed.
This condition determines these fields in a unique manner. 
The author proves that, as a consequence of this principle,  the vector
fields ${\bf E}_j$ and ${\bf B}_j$ must not depend on the velocities of the material points 
and must verify both \emph{Maxwell equations}
 $$ \rot {\bf E}_j+\frac{\partial{\bf B}_j}{\partial t}=0\,,
 \quad \div{\bf B}_j=0\,,\quad 1\leq j\leq N\,.
 $$ 
The author shows that Maxwell's principle is well verified in a lot of cases: 
the $N$-body problem,  {a} material point in 
the gravity field,  {an} electrically charged particle in an external 
electromagnetic field. In this latter case, the above equations are the 
first two Maxwell's equations (Maxwell-Faraday equation and Maxwell-Thomson equation) 
that must be verified by the external 
electromagnetic field. 
The author deduces from them the well-known formula giving the expression of 
the \emph{Laplace force}, and concludes  {that} this force 
\emph{is not} a relativistic effect since it results from the application of 
Maxwell's principle in the framework of
classical mechanics.
\par\smallskip

By contrast, the creation of  {an} 
electromagnetic field by electric charges in motion,
 mathematically described by the last two 
Maxwell equations 
(Maxwell-Ampère equation and Maxwell-Gauss equation) 
are relativistic effects which do not appear  {in the framework of
classical mechanics}. The author adopts Maxwell's principle as a new \emph{law of the mechanics}.
This allows him the study, in the framework of
classical mechanics, systems more general than those made of
material points. For the systems of material points, Maxwell's principle allows, 
under certain conditions, to define a Lagrangian and to show that the Lagrange form 
is nothing else than the exterior derivative of the Cartan form defined in the 
first chapter, in the study of calculus of variations. We can then 
apply the \emph{principle of least action}, often considered as an essential piece of
analytical mechanics. Moreover, when the Lagrangian is hyper-regular, 
one can associate to it a Hamiltonian and use for the study of motions 
the Lagrangian and Hamiltonian formalisms. 
Without denying the importance of the principle of least action nor 
the usefulness of these formalisms, the author declares  that
these concepts seem to him less fundamental than 
Maxwell's principle. His viewpoint seems to him justified
because the existence of a Lagrangian is ensured only locally, and because 
there exist important systems, such as those
made of particles with spin, 
to which Maxwell's principle applies while they have not a globally defined Lagrangian. 
In the sequel, the author will not use the principle of least action except to present, 
very briefly, the \emph{method of variation of constants} introduced by Lagrange 
in 1809 during his works on the slow variations of the orbital elements of
planets.
\par\smallskip

The Lagrange form projects onto the manifold of motions, and its projection 
is a closed form, which is symplectic since it is automatically non degenerate. 
When the system is isolated, the Galilei group acts on %
the evolution space and  {on} the space of motions.  
 {The moment map} of this action 
 {on the evolution space} is a first integral of the motion. 
The author details its ten components, which can be regrouped in three vectors 
$\bf p$, $\bf l$, $\bf g$ and a scalar $E$. He gives their physical meaning:
 $\bf p$ is the total linear momentum; $\bf l$ is the total angular momentum; 
the equality ${\bf g}=$ Constant conveys the fact that the center of mass moves 
on a straight line at constant velocity; the scalar $E$, defined modulo an 
additive constant, is the total energy. 
\par\smallskip

Using a result due to V.~Bargmann \cite{bargmann},
who proved that the symplectic cohomology of the Galilei group is of dimension $1$, 
the author shows that the class of symplectic cohomology of the action of the 
Galilei group  {on} 
the space of motions of an isolated system made of
material points in interaction can be interpreted as the \emph{total mass}. 
The writing of the equations of motion gives two noteworthy results: the vector 
fields ${\bf B}_j$ necessarily all vanish; 
the total force and the total torque of the 
interaction forces are necessarily zero. 
The first result shows that in classical mechanics,
a system made of material points cannot describe
moving magnets. 
The second one expresses the \emph{principle of equality  of action and reaction}, which 
appears as a consequence of 
Maxwell's principle and of the principle of Galilean relativity. 
\par\smallskip

Next the author presents several properties of the dynamical groups of a dynamical 
system and gives numerous examples. For the Kepler problem (motion of a material 
point in a Coulomb field) he shows how to \emph{regularize} 
the manifold of motions and explains the origin of the exceptional first integral, 
often called  the \emph{Lenz vector} or the \emph{Laplace vector}, which should  
rather be called  {the} \emph{eccentricity vector}. 
Although the author does not say it, it should be stressed that this first 
integral was discovered by the Swiss mathematician Jakob Hermann (1678--1833) \cite{Hermann}.
\par\smallskip

In the paragraph \emph{The Principles of symplectic mechanics}, the author first 
works in the  {framework of non relativistic, classical} 
mechanics. He no longer is limited to systems of material points and
adopts the three following assertions as new axioms of 
mechanics: 
\begin{enumerate}

\item{} The space of the motions of a dynamical system is a  
\emph{connected symplectic manifold}.

\item{} If several dynamical systems evolve independently, the manifold of 
motions of the composite system is the \emph{symplectic direct product}
of the spaces of motions of the component systems.

\item{} If a dynamical system is isolated, its manifold of motions admits 
the \emph{Galilei group} as a dynamical group.
\end{enumerate}
 
It is an \emph{extension} of the principles generally admitted in
classical mechanics, which will allow the author to consider new dynamical systems 
having a physical interest.
\par\smallskip

Since, for an isolated system, the author
 {identifies the \emph{mass} of the system with}
the number $m$ which marks 
the class of cohomology of the 
Galilei group action
on the space of motions, 
he can now consider systems of positive, null or negative mass. 
He keeps, for the more general systems that he will consider, 
the physical interpretation of the components of the Galilean moment which he has 
 {previously} given: the vectors ${\bf l}$ and ${\bf p}$
are, respectively, the \emph{angular momentum} and the \emph{linear momentum} of the system,
the scalar  $E$  {is} its \emph{energy}. As for ${\bf g}$, it allows, 
for a system of non vanishing mass $m$, to define the 
\emph{center of mass}  of the system that the author calls  {the}
\emph{center of gravity} or  {the} \emph{barycenter}: 
it is the point whose
position vector, at each instant $t$, is ${\bf R}=({\bf p}t+{\bf g})/m$. 
Hence this point moves  on a straight line at constant velocity  ${\bf p}/m$. 
The author chooses as fundamental quantities the length  $L$, the time $T$ 
and the action $A$, and indicates the dimensionnal equations of the encountered 
quantities in the mathematical description  of the dynamical system: coordinates 
of an element of the evolution space, Lagrange form, components of the matrices 
of an element of the Galilei group, of an element of its Lie algebra and of an element of its dual. 
\par\smallskip

Next, the author proves two theorems concerning the action of a dynamical system 
on a symplectic manifold, that he might have placed in Chapter II. The first provides, 
by means of the action of a dynamical group $G$ on a connected symplectic manifold $V$ 
and the associated symplectic cocycle, the expressions of a moment and that of 
the associated symplectic cocycle of the restriction of this operation to an 
invariant subgroup $G'$ of $G$. The second applies to a moment $\psi$ of the action 
of a connected Abelian dynamical group $G$ on a connected symplectic manifold $V$ 
when the differential $f$ at the identity element of the symplectic cocycle 
$\theta$ associated to this moment is a non degenerate bilinear 
form on the Lie algebra $\mathcal G$ of $G$. 

The author points out that there exist two Abelian subgroups of dimension $3$ 
of the Galilei group of changes of inertial reference frame
of space-time: the group of translations
of the the  {frame of space, the frame of} time 
remaining unchanged, and the group of  {changes of inertial reference 
frame of space-time in which a first inertial frame is replaced by a second one whose 
relative motion with respect to the first one is a translation at a constant velocity,} 
the time frame remaining once again unchanged. The direct product of these two groups 
is a normal subgroup of the Galilei group, isomorphic to the additive group $\RR^6$. 

Applying the above-mentioned theorems, the author shows that the space of motions 
of an isolated system of mass $m\neq 0$ is the  direct product of two spaces 
of motions: the \emph{space of motions of the center of mass}, 
isomorphic to the space of motions of a material point of mass $m$, and the 
\emph{space of motions around the center  of mass}, in which the center of mass 
remains always at the origin of the  reference frame of space.
This important result, called  {the} \emph{barycentric decomposition} of the 
motions of an isolated system, is well known in classical mechanics. 
The Galilei group is a dynamical group  {both of the space of motions 
of the center of mass and of the space of motions around the center of mass}, 
but acts on the last one only through the
quotient by its normal Abelian subgroup isomorphic to 
$\RR^6$. This quotient is isomorphic to $\SO(3)\times\RR$. Hence
the direct product of the Galilei group and of $\SO(3)\times\RR$ is a dynamical group 
of any isolated system of non vanishing mass. The author gives several examples 
including the one, non classical, of a particle with spin.
\par\smallskip

Next the author tackles the study of the relativistic systems. He calls 
\emph{Lorentz frame} an inertial reference frame of
space-time  in which the time unit has been chosen in such a way that 
the velocity of light is equal to $1$. He presents the essential mathematical 
concepts used in the special theory of relativity~: 
the Minkowski space, the Lorentz group and the Poincar\'e group. 
He gives the expressions of a diffeomorphism of $\SO(3)\times\RR^3$ onto 
the restricted Lorentz group ( {the} connected component of this 
group containing the identity) and of 
matrix representations of the Lie algebras of these groups. 
He points out, without proof, several important properties: the Lie algebra of 
the Poincar\'e group is equal to its derived algebra; the symplectic cohomology 
of the Poincar\'e group, and more generally its cohomology valued in the dual 
of its Lie algebra, are trivial.
\par\smallskip

In special relativity, the passage from
a Lorentz frame to another one is made by the action of an element of the 
restricted Poincar\'e group. The author establishes the formula of change
of Lorentz frames and gives its matrix expression. The firts two axioms of
symplectic mechanics remain unchanged, and in the third one, the Galilei group 
is replaced by the restricted Poincar\'e group. Therefore these axions are

\begin{enumerate}
\item{} The space of motions of a dynamical system is a 
\emph{connected symplectic manifold}.

\item{} If several dynamical systems evolve independently, the manifold of motions 
of the composite system is the \emph{ direct product} of the spaces of motions 
of the component systems.

\item{} If a dynamical system is isolated, its manifold of motions admits the 
\emph{restricted Poincar\'e group} as a dynamical group.

\end{enumerate}

It is always possible to choose the moment  {map} of the action 
of the restricted Poincaré group on the space of motions of an isolated 
relativistic dynamical system so that the cocycle associated to the 
moment  {map} is null. 
This condition determines the moment map in an unique way, 
while in classical mechanics, the moment map
of the action of the Galilei group on the space of motions of an isolated system 
depends on an arbitrary additive constant. Moreover, in relativistic mechanics,  
the barycentric decomposition of
motions of an isolated system no longer exists.
\par\smallskip

The intersection of the Galilei and Poincaré groups, considered as subgroups of 
the group of affine transformations of
space-time, is a dynamical group of dimension $7$ of 
isolated dynamical systems, both classical and relativistic.  
The  {moment map} of its action on  the space of motions of an 
isolated system is composed of the vectors ${\bf l}$, ${\bf p}$ and the scalar $E$, 
that the author interpreted in classical mechanics as being the angular momentum, 
the linear momentum and the energy. For a relativistic system, the author chooses 
to conserve for these quantities the same interpretation as in classical mechanics. 
The choice of a Lorentz frame allows to associate to the couple $({\bf p}, E)$  
a vector $P$ of the Minkowski space-time, called the \emph{$4$-momentum} or 
the \emph{energy-momentum vector} of the system. 
Next the author indicates some  {formulae} useful in geometry of 
oriented Minkowski space, refering for more details to his book \cite{CL}. 
The action of the Poincaré group on the space of motions of an isolated 
relativistic system  reveals a second vector $W$ of the Minkowski space-time, 
called  {the} \emph{polarization vector}, orthogonal
to the energy-momentum vector $P$.
\par\smallskip

In a long section, the author proposes a mechanistic description of elementary
particles. In the framework of relativistic mechanics, an isolated dynamical 
system is said to be \emph{elementary} when the Poincaré group acts transitively 
on the space of its motions. The moment map of its action is then a symplectic
diffeomorphism of this space onto a coadjoint orbit of the Poincaré group. 
For the author, the so defined \emph{elementary systems} are 
mathematical  {models for} 
\emph{elementary particles} of physicists. He uses the
type (timelike, spacelike or lightlike, or in other words isotropic) of the
quadrivectors $P$ and $W$, defined in the previous section, for a classification of
elementary systems. By these means, he obtains a large part of the 
 {physicists'} classification of elementary particles. Below, 
briefly summarized, his resultss are presented.   
\par\smallskip

{\bf Case 1, a particle with spin.} It is when $P$ is timelike and when $W$ 
(which, being orthogonal to $P$, is spacelike) is non-zero. The author proves 
that two real numbers $m$ and $s$, whose expressions are given in terms of 
$P$ and $W$, can be interpreted as the \emph{mass} and the \emph{spin} of the 
particle. These numbers are constant on the space of motions of the system. 
The mass $m$ is non-zero, but can be either positive or negative, while the spin 
$s$ is always strictly positive. For each given pair $(m,s)$, with $m\neq 0$ 
and $s>0$, there exists only one model of particle with mass $m$ and spin $s$. 
Its space of motions is $8$-dimensional. For each motion of the system, there 
exists an affine straight line of the Minkowski space-time, parallel to the 
timelike quadrivector $P$, interpreted as the \emph{trajectory} of the particle. 
By expressing the energy-momentum quadrivector $P$ in any Lorentz reference frame, 
the author observes that the norm of the velocity $v$ of the particle in that 
reference frame is always smaller than $1$ (which is, with the chosen units, 
the norm of the velocity of light) and obtains the famous Einstein's formula
 $$E=\frac{m}{\sqrt{1-\Vert {\bf v}\Vert^2}}\,.
 $$
\par\smallskip

{\bf Case 2, a particle without spin.} It is when $P$ is timelike and $W=0$. 
There still exists a real number $m\neq 0$, expressed in terms of $P$, 
interpreted as the \emph{mass} of the particle, constant on the space of motions 
of the system. For each  given real $m\neq 0$, there exists only one model of 
particle without spin with mass $m$, and its space of motions is $6$-dimensional. 
Still in that case, for each motion  of the system, there exists an affine 
straight line of the Minkowski space-time, parallel to the timelike quadrivector 
$P$, interpreted as the \emph{trajectory} of the particle. In any Lorentz 
reference frame the norm of the velocityof the particle is always smaller than 
$1$, the energy $E$ and the velocity $v$ of the particle are related by the above 
written Einstein's formula.
\par\smallskip

 {\bf Case 3, a massless particle.} It is when both $P$ and $W$ are non-zero
and lightlike. The author defines three real numbers $\eta=\pm 1$, $\chi=\pm 1$ 
and $s>0$, interpreted, respectively, as the \emph{sign of the energy},
the \emph{helicity} and the \emph{spin} of the particle, expressed in terms of
$P$ and $W$, constant on the space of motions of the system. For each triple 
$(\eta,\chi,s)$ of real numbers satisfying $\eta=\pm 1$, $\chi=\pm 1$ and $s>0$, 
there exists only one model of massless particle  with sign of the energy $\eta$, 
helicity $\chi$ and spin $s$, and its manifold of motions is $6$-dimensional.
For each motion of the system, there exists a three-dimensional affine subspace 
of the Minkowski space-time which, interpreted in any Lorentz reference frame, 
can be described as a two-dimensional spacelike plane moving at the velocity of light,
called the \emph{wavefront} of the particle. There is no more an affine straight 
line of the Minkowski space-time which can be considered as the trajectory of the
particle. More precisely, for each Lorentz reference frame, there is such 
an affine straight line, which is lightlike and contained in the wavefront of the particle.
However, this straight line depends on the chosen Lorentz reference frame and sweeps 
the whole wavefront when the considered Lorentz reference frame takes all the 
possible values.  
\par\smallskip

The author briefly indicates the existence of other elementary systems, 
for example \emph{tachyons}, which do not correspond to known elementary particles.
Then he looks at \emph{non relativistic elementary particles}, 
beginning by particles without spin. He obtains a mathematical model of such particles
by means of a suitable change of variables in which the velocity of light $c$ appears, and
then by letting  $c\mapsto +\infty$. The model so obtained corresponds to material
points considered in classical non-relativistic mechanics. The same procedure starting 
with relativistic particles with spin leads to a model of non-relativistic material
point with spin, interpreted as a \emph{proper angular momentum}. He briefly indicates another
way in which mathematical models of non-relativitic particles could be obtained,
using as spaces of motions orbits of affine actions of the Galilei group on the dual of its Lie algebra,
which may involve a symplectic cocycle. By this means he obtains models of non-relativistic
massless particles moving at an infinite velocity, each model being characterized by three real numbers
$\chi=\pm 1$, $s>0$ and $k>0$, called, respectively, the \emph{helicity}, the \emph{spin}
and the \emph{color} of the considered particle.      
\par\smallskip

At the end of this section, the author explains that the theory 
of general relativity argues in favour of physical elementary dynamical systems whose
space of motions admits the full Poincaré group $G'$ as a dynamical group. 
Such a system's space of motions may  have several
connected components. The full Poincaré group $G'$ has four connected components, and
all its elements are obtained by composition of an element of the restricted
Poincaré group $G$ (the connected component of the neutral element) with elements 
of two discrete subgroups, each with two elements: the group of \emph{space inversions} 
(exchange of right and left), and the group of  
\emph{time reversals} (exchange of past and future). The author fully discusses 
geometric properties of the moment map of a Hamiltonian action of $G'$ on a
(maybe non-connected) symplectic manifold, as well as geometric  
properties of its coadjoint orbits, and presents their consequences for physical
isolated elementary systems whose space of motions admits $G'$ as a dynamical group. 
He considers elementary  particles first with a non-zero mass, then with zero mass.
Known massless particles (photons and neutrinos)  exist with two opposite helicities, 
and the author considers this fact as an argument in favour of the admission of 
the group of space inversions as an invariance group of mechanics. 
\par\smallskip

Chapter III ends with a study of particles dynamics first in the framework of 
classical mechanics in a fixed inertial reference frame.  Starting from the dynamics
of a free material point, the author explains how to describe the dynamics 
of an electrically charged particle submitted to an electric fild $\bf E$ and a magnetic field
$\bf B$. He uses a modification of the Lagrange form to account for 
the effets of $\bf E$, $\bf B$ and the electric charge of the particle. 
Using the fact that $\bf E$ and $\bf B$ satisfy the
Maxwell equations, he proves that this modified Lagrange form still 
is closed, and that its kernel, of dimension $1$, still determines a foliation 
in curves of the  evolution space, and derives the equations of motion.
The same procedure, in which the symplectic form
of the space of free motions of a relativistic particle is used (instead of that
of the space of free motions of a classical particle), leads to the equations of
motion of a relativistic electrically charged particle submitted to an 
electromagnetic field. In these equations, in good agreement with experimental results,
appears the relativistic correction of the linear momentum well known by physicists.
\par\smallskip

For a particle with spin, either classical or relativistic, submitted to 
an electromagnetic field, the author explains that experiments show that still 
another term must be added to the Lagrange form. This term is the product
by a constant $\mu$ (later interpreted as the module of the magnetic 
moment of the particle) of the exterior derivative of a $1$-form $\varpi$
which, at the non-relativistic limit, is the product of $\d t$ by a function. 
The equations of motion so obtained in the non-relativistic approximation
are in good agreement with the Stern and Gerlach experiments, as well as with 
the precession of spin and the magnetic resonance phenomena. For a relativistic 
particle these equations become very complicated. When the electromagnetic field 
is constant, they are used for the measurement of the 
\emph{anomalous magnetic moment} of particles.  
\par\smallskip

The author now considers $n$ non-interacting particles, which may be either free, 
or subjected to a field. Let $U_i$ be the space of motions  of the system made 
by the $i$-th particle alone. The space of motions
$U$ of the system is the direct product of the spaces of motions $U_i$, 
$1\leq i\leq n$. The author describes two different ways in which an evolution space
can be built for the system, which lead to two different evolution spaces, respectively called
the \emph{synchronous evolution space} and the \emph{asynchronous evolution space}.
The synchronous evolution space is obtained by adding to the space of motions $U$
one dimension for the time, and by taking the initial conditions of all particles 
at the same time. Its dimension is $\dim U+1$. The asynchronous evolution space,
built by taking the initial conditions of the particles at different times, 
is of dimension $\dim U+n$, and involves $n$ different times, one for each particle. 
Both the synchronous and the asynchronous evolution spaces are presymplectic manifolds
which project onto the space of motions $U$. The synchronous evolution space
should be used only for non-relativistic particles since it involves a notion of simultaneity, 
which in relativistic physics depends on the choice of a Lorentz reference frame.
When all the particles of the system are identical and submitted to the same field,
the spaces of motions $U_i$, for all $i\in\{1,\ldots,n\}$, are equal to $U_1$ and
one could think that the space of motions of the system is $(U_1)^n$. 
However, experiments have shown that motions which only differ by the labelling 
of particles should be identified. The true space of motions of the system
is the set of equivalence classes of $n$-uples $(x_1,\ldots, x_n)\in (U_1)^n$ which
satisfy $x_i\neq x_j$ for $i\neq j$, $1\leq i,j\leq n$, two $n$-uples
$(x_1,\ldots,x_n)$ and $(x'_1,\ldots,x'_n)$ being equivalent if there exists a permutation $\sigma$
of $\{1,\ldots,n\}$ such that $x'_i=x_{\sigma(i)}$ for all $i$, $1\leq i\leq n$.           
 The author indicates the expression of the Lagrange form and remarks that the use of a set
of equivalence classes as a set of motions of the system is a consequence of the 
indiscernability of particles which does not involve quantum mechanics.
\par\smallskip

A \lq\lq classical\rq\rq\ method for obtaining the equations of motion of 
a system of interacting particles uses the synchronous evolution space of the system
when the particles do not interact, with a Lagrange form modified by addition of
a suitable term involving an \emph{interaction potential}. This method works 
successfully for celestial mechanics. The author uses it for a system of 
non-relativistic particles with spin, 
indicates the equations of motion and, when the system is isolated, writes the
formulae which express the constancy of the Galilean moment map. The equations 
so obtained take into account the electrostatic and magnetostatic forces exerted 
by each moving particle on the other particles. However, they do not account for other
small, but measurable relativistic effects, such as the Laplace force. To account for 
these effects, the author considers the use of the asynchronous evolution space. 
He proves that the use of that space would prohibit the  existence of interactions 
between the particles, and concludes that a way to get around this difficulty 
could be to abandon the idea of localized particles, the distinction between particles
and fields being a non-relativistic approximation.
\par\smallskip

The \emph{scattering theory} is an approximate mathematical description of
a system of interacting particles in which it is assumed that when each particle
is far enough from all other particles, the interaction forces can be neglected.
The author presents a mathematically rigorous version of this theory. He assumes that
unscattered and scattered dynamical systems share the same evolution space with two
different Lagrange forms, equal of an open subset $\Omega$ of the evolution space.
The complement of $\Omega$, on which the Lagrange forms of the scattered and unscattered
systems are not equal, is called the \emph{scattering source}. A scattered motion 
is said to be \emph{constrained} when it is wholly contained in the scattering 
source. The author looks at unconstrained scattered motions contained in the 
scattering source only for a bounded time interval. He states, without a 
complete proof, a theorem according to which such a motion coincides with two 
different unscattered motions, one before it enters 
the scattering source and another one after it has finally left the scattering source.
The so established correspondence between two unscattered motions is a symplectomorphism
between two open subsets of the set of motions of the unscattered system.
The author defines the \emph{symmetry group} of the scattering source and proves 
that it is a dynamical group for both systems, scattered as well as unscattered.
For \emph{bounded}, \emph{static} or \emph{conservative} scattering sources, some properties
of this group can be deduced. The author briefly presents the dynamical system made 
by photons in a refracting telescope. Due to a term containing a length 
(later interpreted as the wavelenth of the light) in the formulae he obtains,
the image of a distant star can never be a point. When this term is neglected, 
one obtains the \emph{geometrical optics} approximation. Using  the 
scattering theory, the author discusses, in relativitic physics, the reflection 
of light on a moving mirror and obtains formulae for the \emph{Doppler effect due to
reflection}.    
\par\smallskip

Collisions of relativistic free particles with non-zero masses are finally 
discussed by the author, with the help of the scattering theory. Though no valid model of 
relativistic interacting particles is available, the author obtains some properties 
of the symplectomorphism which relates the free motions of the particles 
before and after their collision. He proves that this symplectomorphism commutes
with the action of elements of the Poincaré group, which implies that the total momentum
and the total energy of the system of particles are conserved by collisions.
As a conclusion, he states that knowing this symplectomorphism could allow
to study the \emph{constrained motions} by the technique of analytic continuation.  

\subsubsection*{Chapter IV, Statistical mechanics.} It contains two sections. 
The first one, of about 50 pages, is essentially mathematical. The second one,
of about 35 pages, presents the principles of statistical mechanics.
\par\smallskip

The first section of this chapter begins by introducing various concepts related
to smooth manifolds, topological spaces, ordered vector spaces 
(called \emph{Riesz spaces}), normed vector spaces, specially Banach and 
Hilbert spaces. A very condensed course (about 30 pages) in Measure Theory
and Integration follows, with some notions in Probability. The author uses 
the presentation, privileged by the Boubaki school, in which a measure on a 
smooth Hausdorff manifold $V$ is an element of the topological dual vector space 
of the space of continuous, compactly supported functions on $V$, instead of defining
first measurable parts of $V$, and then a measure as a denombrably additive function
defined on the set of measurable parts. Then the author defines 
\emph{probability measures}. A measure $\lambda$ on $V$ is said to be \emph{defined
by an everywhere positive continuous field of densities} on $V$, when its 
expression in every chart is the product, by an everywhere strictly positive
continuous function,  of the Lebesgue measure. The measure $\lambda$ being fixed, 
the author considers a probability measure whose density with respect to $\lambda$
is a continuous function $\rho\geq 0$. He defines the $\lambda$-\emph{entropy}
of this probability measure by setting
 $$s_\lambda(\rho)=\begin{cases}
           \displaystyle\int_V-\rho(x)\log\bigl(\rho(x)\bigr)\lambda(\d x)&
            \text{if this integral converges,}\\
             \displaystyle-\infty&\text{if the above integral does not converge\,.}
              \end{cases}
 $$            
By convention, at points $x\in V$ where $\rho(x)=0$, the value taken by the function  
$x\mapsto -\rho(x)\log\bigl(\rho(x)\bigr)$ is $0$.
\par\smallskip

A continuous map $\Psi$, defined on $V$, with values in a finite-dimensional 
vector space $E$ being given, the author calls
\emph{generalized Gibbs probability law} any completely continuous probability
law for which the map $\Psi$ is integrable, whose density $\rho$ with respect 
to $\lambda$ is expressed as
 $$\rho(x)=\exp\Bigl(-\bigl(z+\langle Z,\Psi(x)\rangle\bigr)\Bigr)\,,\quad\hbox{with}\ x\in V\,.
 $$
In the above equality, $Z$ is an element of the dual vector space $E^*$ of $E$, 
which must be such that the integrals below, appearing on the right hand sides
of the equalities defining $I_0$ and $I_1$, are convergent. The real $z$ must be chosen
in such a way that the integral of $\rho$ on $V$ with respect to the measure $\lambda$
is $1$.
\par\smallskip

The  \emph{normal distribution} on $\RR^n$ is a generalized Gibbs probability law 
for a suitable choice of the map $\Psi$.
\par\smallskip

The mean value of $\Psi$, for the probability law of density $\rho$ with respect to
$\lambda$, being denoted by $M$, and setting
 $$ 
  I_0=\int_V\exp\Bigl(-\langle Z,\Psi(x)\rangle\Bigr)\lambda(\d x)\,,\quad
  I_1=\int_V\Psi(x)\exp\Bigl(-\langle Z,\Psi(x)\rangle\Bigr)\lambda(\d x)\,, 
 $$
the author can write
 $$z=\log(I_0)\,,\quad M  =\int_V\Psi(x)\rho(x)\lambda(\d x)=\frac{I_1}{I_0}\,.
 $$
The $\lambda$-entropy of the generalized Gibbs law of density 
$\rho$ with respect to $\lambda$ is
 $$s(\rho)=z+\langle Z,M\rangle\,.
 $$
The author proves that, on the set of completely continuous probability laws for which 
the mean value of $\Psi$ is $M$,  the $\lambda$-entropy functional has a strict maximum
at the generalized Gibbs probability law of density $\rho$ with respect to $\lambda$.
Moreover, when the set of values of $\Psi$ is not contained in an affine subspace of
 $E$ of dimension strictly smaller than $\dim E$, he proves that the map $Z\mapsto M$
is injective. He derives conditions for the differentiability  with respect to 
$Z$, under the sign $\int$, of the integrals $I_0$ and $I_1$. He will later improve
these results in his paper \cite{Souriau1974}, 
published a few years after his book \cite{SSD}.   
\par\smallskip

The author then assumes that the  Hausdorff manifold $V$ is endowed with a symplectic form 
and that a connected Lie group $G$ acts on it by a Hamiltonian action. For $\lambda$,
he chooses the Liouville measure, and for the map $\Psi$ a moment map of the action of $G$.
The vector space $E$ is therefore now the dual vector space ${\mathfrak g}^*$ 
of the Lie algebra $\mathfrak g$ of $G$, and its dual $E^*$ is ${\mathfrak g}$.
He denotes by $\Omega$ the largest open subset of $\mathfrak g$ on which
the integrals $I_0$ and $I_1$ are convergent and define differentiable functions of 
the variable $Z$ whose differentials are continuous and can be obtained by 
differentiation under the $\int$ sign. A generalized Gibbs probability law can then be 
associated to each $Z\in \Omega$. The author indicates the corresponding expressions of 
$z=\log(I_0)$, of the entropy $s$ and of the mean value $M$ of the moment map, 
considered as functions of the variable $Z\in{\mathfrak g}$. He proves that when
$G$ acts effectively on $V$, the map $Z\mapsto M$ is injective and open, therefore is
a diffeomorphism of $\Omega$ onto an open subset $\Omega^*$ of ${\mathfrak g}^*$. 
By using the inverse map, $z$, $s$ and $Z$ can be considered as functions of the variable 
$M$, which spans the open subset $\Omega^*$ of ${\mathfrak g}^*$. The author 
proves that these functions can be differentiated and indicates the expressions 
of their differentials. Adding a constant to the moment map does not change 
the generalized Gibbs laws on the manifold $V$. Their set is called by the author
the \emph{Gibbs set} of the considered Lie group action. Moreover the author proves 
that the open subset $\Omega$ of $\mathfrak g$ is a union of orbits of the 
adjoint action of $G$ and is endowed with a definite negative Euclidean metric.      
\par\smallskip

Results obtained by the author in the first section of Chapter IV are, in 
the second section, applied to dynamical systems encountered in physics. 
A \emph{statistical state} of such a system is a probability law on the space 
of its motions. The author explains that in the \emph{kinetic 
theory of gases}, a gas in a container at rest in a Galilean reference frame
is considered as an assembly of very small material particles whose motions are 
governed by the laws of classical mechanics, interacting by instantaneous 
elastic collisions between themselves and with the walls of the container. 
With these assumptions, the gas can be modelled by a conservative Hamiltonian 
system. The one-dimensional group of time translations acts on the symplectic 
manifold of motions of the system by a Hamiltonian action, with the energy as 
a moment map. The author explains that the entropy of the system increases with time,
so assuming that the \emph{natural equilibria} of the gas are elements of the Gibbs set 
of the group of time translations is a very reasonable assumption. 
Each Gibbs state is determined by an element $Z$ of the one-dimensional
Lie algebra of this group, which is a way of measuring the \emph{gas temperature}.   
Each Gibbs state is unaffected by the action of the group of time translations, since the
adjoint action of this group on its Lie algebra is trivial. The author indicates the
dimension equations of quantities which determine a Gibbs state, or appear 
in its description. He explains that when a compound system is the union of several subsystems,
the tensor product of natural equilibria of the subsystems is a natural equilibrium of
the compound system only when the element $Z$ of the Lie algebra of the group of time translations
which determine the natural equilibria is the same for all subsystems, in other words
only when the temperatures of all the subsystems are equal. If this condition is 
not satisfied, the compound system is in a state of natural equilibrium only when 
the subsystems cannot exchange energy between them. As soon exchanges of 
energy can occur, even when they are very tiny, the compound system is no more 
in a natural state of equilibrium.
\par\smallskip

For the determination of the Gibbs states of a gas, one has to calculate the integral 
$I_0$. The \emph{ideal gas approximation} is when, for this calculation, the only 
motions taken into account are those for which, at the considered time, no 
collisions occur between particles or between a particle and the walls of the container, 
and when the total volume of the particles is considered as negligible
in comparison with the volume of the container. The author successively considers
monatomic and polyatomic ideal gases and derives, for a Gibbs state, the 
probability distribution of velocities of the particles, which was determined by 
Maxwell in 1860. He explains the principle of an \emph{ideal gas thermometer} used
for precise measurements of the temperature. He is led to the formula $Z=1/(kT)$, which
expresses the absolute temperature $T$ in terms of the element $Z$ of the Lie algebra
of the group of time translations, and of \emph{Boltzmann's constant} $k$. He obtains
the Mariotte, Gay-Lussac and Avogadro laws for perfect gases well known by physicists.
He proves that for an ideal gas thermometer, the temperature is a random variable whose
probability law is very concentrated around its mean value, and converges
(for the weak topology)  towards the Dirac measure when the number of particles 
tends to infinity.  
\par\smallskip

For a conservative system, the entropy $s$ of a Gibbs state is a smooth function 
of the energy $E$ whose derivative is the element $Z$ of the Lie algebra which 
determines the Gibbs state. Units being chosen, $Z$ becomes a real 
number. Experimental results show that this number is always strictly positive. 
The author can therefore consider the energy $E$ of
a Gibbs state as a function of its entropy $s$ and of other variables which 
describe the system. For example, when the considered system is a gas, 
these additional variables are the volume of the container, the number of particles, 
etc. The author assumes that these variables can be globally described 
by an element $u$ in a smooth Hausdorff manifold. Infinitesimal variations of a 
Gibbs state are therefore described by the differential of $E$,
 $$\d E= \frac{\d s}{Z}+\varpi\d u=T\d S+\varpi\d u\,,
 $$   
where $\varpi$ denotes the partial differential of $E$ with respect to $u$, and
$S=ks$ ($k$ being Boltzmann's constant). The author calls $S$ and $u$ the
\emph{position variables}, $T$ and $\varpi$ the associated \emph{tension variables}. 
For example, when the considered system is a gas in a container whose volume may vary,
the above expression becomes
 $$\d E= T\d S-p\d V\,,
 $$
where $p$ is the \emph{pressure}.
\par\smallskip

The author then explains that the infinitesimal variation of the energy $E$ is the sum of 
the infinitesimal variations of two quantities $Q$ and $W$, called respectively
the \emph{heat} and the \emph{work}, which are not functions, but rather
\emph{action integrals} in the sense of calculus of variations. Variations of
heat and work depend on the path followed while the adjustable variables
which describe the state of the system changed. Different types of evolution of 
the state of a system should be separately considered, for example 
\emph{adiabatic evolution}, \emph{isothermal evolution}, etc.
\par\smallskip

The author defines the \emph{heat capacities} of a gas at \emph{constant volume}
and at \emph{constant pressure}, and indicates the expressions of the 
thermodynamic functions of an \emph{ideal fluid}. These expressions depend on the
chosen model for the particles. A table indicates their expressions and values for 
various models (material point, particle with spin, etc). Agreement 
with experimental results is good for a monoatomic gases, not so good 
for polyatomic gases and for solids at low temperatures. 
\par\smallskip

Since the group of time translations is a subgroup, but not a \emph{normal}
subgroup, of the Galilei group, a dynamical system conservative in some inertial
reference frame is not conservative in a different inertial reference frame.
This important remark leads the author to introduce the new concept of 
\emph{covariant statistical mechanics} by proposing the following principle:
\par\smallskip

\emph{When a dynamical system is invariant by the action of some Lie subgroup 
$G'$ of the Galilei group, its natural equilibria are the elements of the Gibbs 
set of the action of $G'$.}
\par\smallskip

Each natural equilibrium of such a system is determined by an element $Z$ of 
the Lie algebra ${\mathfrak g}'$ of $G'$, which is a Lie subalgebra of the Lie 
algebra of the Galilei group. The author observes that $Z$ generalizes the inverse 
of the temperature, and discusses its physical interpretation. For this purpose he
chooses an inertial reference frame ${\mathcal R}_0$ of the Galilei space-time and 
considers the reference frame ${\mathcal R}=\exp(tZ){\mathcal R}_0$, which 
generally is not inertial. He looks at the indications given, in the reference 
frame $\mathcal R$, by an ideal gas thermometer in equilibrium with the considered 
system. It amounts to observe the system in a moving frame. The author proves that
it is as if $\mathcal R$ were an inertial frame and if the particles of gas were 
submitted to additional forces (inertial forces, \emph{i.e.}, centrifigal and 
Coriolis forces). By this means he obtains general formulae for the probability 
density of the Gibbs state associated to $Z$. He then discusses in greater detail
several examples: the wind (the relative motion of $\mathcal R$ with respect to 
${\mathcal R}_0$ is a translation at constant velocity), accelerating rocket
(this relative motion is a uniformly accelerated translation), gas in a centrifuge 
(this relative motion is now a rotation around an axis at a constant angular velocity).
In this last example he considers also a system made by particles with spin, and finds that the
most probable orientation of the particle's spin is parallel to the rotation axis.
\par\smallskip  

One may wish to apply the above principle for a system invariant 
by the whole Galilei group. However, the subset $\Omega$ of the Lie algebra of 
the Galilei group made by elements which determine a Gibbs state is then empty.
Looking at the motions of the system around its center of mass, the author is led 
to consider the just discussed equilibria in a rotating frame. For the author,
the rotation of celestial bodies observed in Astronomy confirms the validity of 
his principle of covariant statistical mechanics.   
\par\smallskip

For a relativistic dynamical system, the Galilei group must 
be replaced with the Poincaré group. As above, the Gibbs state of the action of 
the whole Poincaré group is empty, and the author is led to consider systems invariant 
by the action of a subgroup $G'$ of the Poincaré group. He presents in more detail 
several examples: an ideal gas in a container at rest in an inertial frame
(in that example Maxwell's distribution law for the velocities of particles is
slightly modified); relativistic wind; statistical equilibria of photons. In this 
last example, the number of photons cannot be fixed. The author explains how 
such a system can be described and even takes into account the fact that photons
can have two opposite circular polarizations. However, his formula is not 
in agreement with the one obtained by Planck for the black-body radiation, which  
is in very good agreement with experimental observations. The author concludes 
that this last example must be dealt with in the framework of quantum mechanics.
\par\smallskip

\subsubsection*{Chapter V, Geometric quantization.} 
It contains two long sections, each of around 40 pages. Before 
describing their contents, let us recall
that a \emph{contact form} on a smooth Hausdorff manifold $Y$ is a 
differential $1$-form $\omega$, which nowhere vanishes on $Y$, 
whose exterior derivative $\d \omega$ is non degenerate on $\ker\omega$. 
The existence of a contact form on the manifold $Y$ has several important
consequences: $Y$ must be odd-dimensional, the $2$-form $\d \omega$ is a
\emph{presymplectic form} and there exists on $Y$ a unique vector field
$R_Y$, called the \emph{Reeb vector field} (in honour of the French mathematician
Georges Reeb), determined by the two equalities  $\mathi(R_Y)(\d\omega)=0$ and 
$\mathi(R_Y)\omega=1$.  
\par\smallskip

The author defines a \emph{quantum manifold} as a smooth manifold $Y$ 
endowed with a contact form $\omega$ such that all integral curves of the
Reeb vector field $R_Y$ are $2\pi$-periodic. These curves are then the orbits 
of an operation on $Y$ of the one-dimensional torus. The set of these curves,
in other words the quotient of $Y$ by this operation, is a symplectic manifold
$(U,\sigma_U)$, called by the author the \emph{basis} of the quantum manifold $Y$.
The quantum manifold projects on $U$, and the pull-back by the projection map
of the symplectic form $\sigma_U$ is the presymplectic form $\d \omega$.  
\par\smallskip

A \emph{quantization} of a given Hausdorff symplectic manifold $(U,\sigma_U)$
is defined by the author as the construction of a quantum manifold whose basis is
$(U,\sigma_U)$. A symplectic manifold is said to be \emph{quantizable} when its 
quantization is possible. The author indicates, without proof, a necessary and sufficient
condition in which a  given symplectic manifold is quantizable: the cohomology class
of its symplectic form must be $2n\pi$, with $n$ an integer. It is satisfied
when $\sigma_U$ is the exterior derivative of a $1$-form. The manifolds of motions
of many mechanical systems, for example harmonic oscillators, Kepler's problem,
and the $N$-body problem of celestial mechanics, are therefore quantizable.
The author proves that the manifold of motions of a non-relativistic particle with spin 
is quantizable if and only if the spin of the particle is integer or half integer,
when expressed with $\hbar$ as unit. The author indicates several examples of
quantizable symplectic manifolds, together with the full description of the 
corresponding quantum manifolds: two-dimensional spheres \textcolor{black}{of integer or half-integer radii,}
spaces of motions of a relativistic particle, first with a non-zero mass, then 
with zero mass, with spin $1/2$ or $1$.      
\par\smallskip

Isomorphisms of quantum manifolds are called by the author \emph{quantomorphisms}.
Any quantomorphism between two quantum manifolds projects onto a symplectomorphism between their
bases. Two quantizations $Y$ and $Y'$ of the same symplectic manifold are said to be \emph{equivalent}
when there exists between them a quantomorphism which projects onto the identity map of their common basis.
A symplectic manifold is said to be \emph{monoquantizable} when all its quantizations are equivalent.
\par\smallskip
  
A group $\Gamma$ of quantomorphisms of a quantum manifold 
$Y$ projects onto a group of symplectomorphisms of its basis $U$, and its
projection  is a group homomorphism. Conversely, a group $G$ of symplectomorphisms of 
the basis $U$ is said to be \emph{liftable} if there exists a group $\Gamma$
of quantomorphisms of $Y$ which projects onto it. When in addition the projection
of $\Gamma$ onto $G$ is a group isomorphism, $\Gamma$ is said to be an 
\emph{isomorphic lift} of $G$. The author proves that the set of isomorphic lifts
of $G$ is in one to one correspondence with the set of its characters
(group homomorphisms of $G$ into the one-dimensional torus).
\par\smallskip
A theorem due to the author states that a simply connected 
quantizable symplectic maniold is monoquantizable. Another theorem explains
how to  quantize a symplectic manifold $U$ when a quantization of a covering 
manifold $U'$ of $U$ is known. Conversely, when a quantization $Y$ of $U$ 
is known and when the group of symplectomorphisms of $Y$ determined by 
the covering manifold $U'$ isomorphically lifts onto a group $\Gamma$ of 
quantomorphisms of $Y$, $\Gamma$ is a discrete group and can be used 
to build a covering manifold $Y'$ of $Y$, which quantizes $U'$. Using these 
two theorems, it can be proven that there exist as many non equivalent
quantizations of a quantizable connected symplectic manifold as 
its homotopy group has distinct characters. Therefore, being simply connected,
the two-dimensional sphere, any symplectic vector space, the space of motions of 
a free particle, (with or without spin, non-relativistic or relativistic),
when quantizable, are monoquantizable. The author then discusses the quantizability
of the manifold of motions of a system of $N$ particles without interactions. If each of these 
particles can be distinguished from the others and has an integer or half-integer spin, 
the space of motions of the system is quantizable and simply connected, therefore monoquantizable.
But if each particle cannot be distinguished from the others, the space of 
motions of the system has exactly two non-equivalent quantizations, 
corresponding to the two distinct characters of the group of permutations of a set of
$N$ elements. These two quantifications are the Bose-Einstein and Fermi-Dirac 
quantizations, well known to physicists.
\par\smallskip

A smooth vector field, defined on a quantum manifold $Y$, whose flow acts on
$Y$ by quantomorphisms, is called by the author an 
\emph{infinitesimal quantomorphism}. The Lie derivative of the contact form $\omega$
defined on $Y$ with respect to an infinitesimal quantomorphism vanishes.
The author proves that each infinitesimal quantomorphism is associated to a smooth function
defined on the basis $U$ of $Y$. For example, the Reeb vector field is an
infinitesimal quantomorphism associated to the constant function whose 
value at any point in $U$ is $1$.  The Lie bracket of two infinitesimal 
quantomorphisms is an infinitesimal quantomorphism, whose associated function is the 
Poisson bracket of  the functions associated to these two infinitesimal quantomorphisms.
\par\smallskip

The author then discusses the \emph{quantization} of a dynamical group
of a quantizable symplectic manifold $U$, with the quantum manifold $Y$ 
as quantization. When a Lie group acts on $Y$ by quantomorphisms, it acts also on 
the basis $U$ by symplectomorphisms, the projection of $Y$ onto $U$ being 
equivariant with respect to these actions. Therefore $G$ is a dynamical group of $U$,
and is said to be a \emph{quantizable} dynamical group. A quantizable 
dynamical group of $U$ is always liftable. The author proves that when
a dynamical group of $U$ is quantizable, its symplectic cohomology is zero. 
He gives examples which prove that this necessary condition is not sufficient.
\par\smallskip

A dynamical group of $U$ which is liftable, but not quantizable, may have
an extension which still is a dynamical group of $U$ and is quantizable. 
As a first example, the author considers a symplectic vector space $E$
which acts on itself by translations. It is a dynamical group of $E$. The product
$Y=E\times T$, where $T$ is the one-dimensional torus, is a quantization of $E$,
which happens to be a Lie group, called the \emph{Weyl group}, extension of
the Abelian group $E$. It acts on itself by translations, which are quantomorphisms.
The Weyl group is therefore a quantizable extension of the Abelian group $E$, 
which itself is not quantizable since its symplectic cohomology is not zero.  
As a second example, the author considers the group $\SO(3)$ acting on 
a two-dimensionl sphere $S_2$ centered on the origin of $\RR^3$, 
endowed with its element of area form as symplectic form. This symplectic 
manifold is quantizable if its total area is an integer or an half-integer.
In this latter case, the author proves that the dynamical group $\SO(3)$, which is liftable
since $S_2$ is simply connected, is not quantizable, though its symplectic 
cohomomogy is zero. He proves also that the group $\SU(2)$ is a quantizable 
extension of $\SO(3)$.  A similar situation occurs for the restricted
Poincaré group, which is a dynamical group of the manifold of motions 
of a relativistic particle. This group is quantizable if and only if the
particle's spin is integer. If the particle's spin is half-integer, 
there exists a quantizable extension of the restricted Poincaré group descibed by
the author in terms of \emph{Dirac's spinors}.   
\par\smallskip

The author denotes by ${\mathcal H}(Y)$ the vector space of smooth, 
complex-valued and compactly supported functions defined on a quantum manifold $Y$,
which are equivariant with respect to the one-dimensional torus $T$ actions
on $Y$ by the flow of the Reeb vector field and on $\CC$ by multiplication 
(the torus $T$ being identified with the set of complex numbers of modulus $1$).
After indicating the definitions and results in topology and
functional analysis he is going to use, he proves that 
${\mathcal H}(Y)$ is a pre-Hilbert space. He defines the \emph{Hilbert space of
the quantum manifold $Y$} as the completion of ${\mathcal H}(Y)$. He defines also
$C^*$-algebras, indicates some of their properties and proves that the set of
bounded linear endomorphisms of a Hilbert space is a $C^*$-algebra. He 
very shortly presents many concepts about operators on a Hilbert space:
self-adjoint, normal, Hermitian, unitary operators, etc.  For a given
quantum manifold $Y$ of basis $U$, the author then explains how to associate
a Hermitian operator  $\widehat{\mathstrut u}$ on  ${\mathcal H}(Y)$ 
to any smooth function $u$ defined on the basis $U$. He proves that the 
map $u\mapsto \widehat{\mathstrut u}$ is linear and injective, that 
for $u=1$, $\widehat{\mathstrut u}$ is the identity
of ${\mathcal H}(Y)$, and that for any pair $(u,u')$ of smooth functions 
on $U$, 
$$\widehat{\mathstrut u}\circ\widehat{\mathstrut u'}-
   \widehat{\mathstrut u'}\circ\widehat{\mathstrut u}
   =-i\widehat{\mathstrut\{u,u'\}}\,.
 $$
He illustrates this important result in the case when $Y$ is the quantization of a
symplectic vector space, the smooth functions on this vector space being 
canonical coordinates. The operators associated to smooth functions in this way
often are not bounded, which makes their study and use very difficult.
He proposes to build directly quantomorphisms acting on $Y$, 
without taking the detour of infinitesimal quantomorphisms and 
exponentiation. He proves a theorem which asserts the existence of an 
injective homomorphism of the group of quantomorphisms of $Y$
into the group of unitary operators on the Hilbert space of $Y$.
\par\smallskip

The second section of Chapter V begins with the statement and explanation of
the \emph{Correspondence Principle}, well known by physicist. Quantum mechanics 
is not an autonomous theory: it cannot be formulated without using classical mechanics.
According to the Correspondence Principle, for each 
physical phenomenon described in the framework of quantum mechanics, there 
should exist a corresponding \lq\lq classical approximation\rq\rq\  described
in the framework of classical mechanics. For the author, this principle 
is an argument in favour of the extension of classical mechanics he proposed 
in Chapter II, the \emph{Maxwell Principle}. Without this extension,  the 
use of classical mechanics would be limited to systems made of material 
points, excluding many phenomena encountered in physics such as particles with spin.
\par\smallskip  

The \emph{quantization} of a classical mechanical system is the construction
of the quantum mechanical system of which this classical system is an approximation.
In this section the author examines the possible application of the mathematical theory
of Geometric Quantization developed in the previous section, 
to the quantization of real systems encountered in physics. His 
initial assumption is: the space $U$ of classical motions of the system is 
a \emph{quantizable symplectic manifold}.  The first consequence of this assumption
is: expressed with the constant $\hbar$ as unit, the value of the spin of 
particles must be either an integer or a half-integer. This consequence 
is in perfect agreement with experimental results. 
\par\smallskip

Let the quantum manifold $Y$ be a quantization of the classical manifold 
of motions $U$. A \emph{state vector} of $Y$ is any element $\Psi$ with
norm $1$ of the pre-Hilbert space ${\mathcal H}(Y)$. The author
explains its probabilistic interpretation:
for each $\xi\in Y$, the real non-negative number $\vert\Psi(\xi)\vert^2=\overline{\Psi(\xi)}\Psi(\xi)$
depends only on the projection $x$ of $\xi$ on the manifold of motions $U$;
the function so defined on $U$ is the probability density, 
with respect to the Liouville measure, of a probability law on $U$ which is the
\emph{statistical state} of the system (in the sense defined in Chapter IV) corresponding
to the state vector $\Psi$. 
\par\smallskip

An \emph{observable} in classical mechanics is a smooth function $u$ defined on the
classical manifold of motions $U$. According to the British scientist
P.A.M.~Dirac, \emph{observables} in quantum mechanics are Hermitian operators.
In agreement with this view, the author has proven in the previous section
that a Hermitian operator $\widehat{\mathstrut u}$ could be associated 
to the classical observable $u$. This operator is the quantum observable 
of which $u$ is the classical approximation. When applied to the phase space of a
classical conservative system at a given time $t$, this idea leads to
Dirac's \emph{quantum equations of motion} and to the \emph{commutation relations}
chosen by the physicists W.~Pauli and W.~Heisenberg as fundamental ingredients of
quantum mechanics.
\par\smallskip

In favorable cases, an isotropic foliation of the classical space of motions
$U$ can be lifted into a \emph{Planck foliation} of the quantum manifold $Y$, 
that means a foliation such that the contact form defined on $Y$ vanishes 
on all vectors tangent to the leaves. When a state vector $\Psi$ is constant
on each leaf of this Planck foliation, this vector state is said by the author
to  \emph{satisfy Planck's condition} (relative to the considered Planck foliation).
The \emph{Planck space} is the set of state vectors which satisfy Planck's condition.  
For a conservative system, using as isotropic foliation of $U$ 
the foliation whose leaves are the integral curves of the Hamiltonian 
vector field whose Hamiltonian is the energy, by writing explicitly Planck's condition,
the author obtains the famous equality $E=h\nu$ which relates the 
energy $E$, the Planck constant $h$ and the frequency $\nu$.
\par\smallskip

The author then assumes that a dynamical group $G$ acts on the classical
manifold of motions $U$, with a moment map, and that there exists 
an Abelian normal subgroup $\widetilde G$ of $G$ whose symplectic cohomology is zero.
He uses as isotropic foliation of $U$ the foliation whose leaves are tangent to
the vector sub-bundle of $TU$ generated by Hamiltonian vector fields whose Hamiltonians
are the components of the moment map of the action of $\widetilde G$. He lifts this foliation into
a Planck foliation of $Y$, writes explicitly the corresponding Planck's condition, and
determines the associated Planck space. For an isolated system, $G$ will be either 
the restricted Poincaré group or the Galilei group, according to whether the system 
is relativistic or non-relativistic. The Abelian subgroup $\widetilde G$ 
will be the group of space-time translations. The author explains how Planck's 
condition leads to the quantum wave equations. For a non-relativistic material point,
the author obtains by this means the \emph{Schrödinger equation}, and for 
a relativistic material point, the \emph{Klein-Gordon equation}. For a 
non-relativistic particle with spin $1/2$, using a $\CC^2$-valued variable 
as state vector, the author proves that both its components
satisfy a Schrödinger equation; it is the description of electrons proposed by Pauli.
For a relativistic particle of spin $1/2$, he obtains \emph{Dirac's equation} and proves that
rather than the Poincaré group itself, it is a quantizable extension of 
this group whih acts on the set of solutions of this equation. For him, the non-quantizability
of the Poincaré group provides a natural explanation of this fact, well known to physicists. 
The author applies the same procedure for a massless particle, first with spin $1/2$, then with spin $1$.  
\par\smallskip

Next the author considers an assembly of particles of a given type in indeterminate number.
He denotes by $U_\Phi$ its manifold of motions and by $U$
the manifold of motions of the system made by a single particle of the considered type.
He calls $U_\Phi$ the \emph{Fock's manifold} and explains that as a set,
it is the set of all finite subsets of $U$. As a manifold, $U_\Phi$ 
is the sum of disjoint parts of different dimensions, each part $U_n$ being 
the manifold of motions of the system made of a fixed number $n$ of particles,
with $0\leq n<+\infty$. For $U_0$, the author takes a singleton considered 
as a $0$-dimensional manifold. As seen in the previous section, for each 
$n\geq 2$, the part $U_n$ of $U_\Phi$ has two non-equivalent quantizations,
one for each distinct character of the group of permutations of a set of $n$ elements. 
The author explains that for some physical considerations, the same character should 
be chosen for all possible values of $n$, therefore he obtains two 
non-equivalent quantizations $Y_\Phi$ of $U_\Phi$. On the pre-Hilbert space 
${\mathcal H}(Y_\Phi)$, he defines \emph{creation} and \emph{annihilation} 
operators, which respectively increase or decrease the number of particles by one unit.
Two different cases must be distinguished, depending on whether the 
particles are \emph{fermions} or \emph{bosons}. 
\par\smallskip

At the end of Chapter V, the author discusses the notion of \emph{average value}
of an observable for a given state vector. Planck's condition appears as a
sufficient condition so that the statistical and quantum mechanical
averages of an observable coincide. Finally the author uses the notion of
\emph{function of positive type} defined on a group, with values in $\CC$,
to enlarge the definitions of  a \emph{quantum state} and of \emph{average value}
of an observable, to include states defined by a \emph{density operator}
encountered in quantum chemistry.

\section{Comments}

The great originality of the book  \emph{Structure des systèmes dynamiques}
clearly appears in the detailed presentation of its contents given above.
In this section, we will first try to identify the most innovative and promising 
ideas which can be found in it. Then we will write a few words about the 
terminology and the notations used by the author, which do not lack originality too.

\subsection*{Remarkabe scientific concepts.} 

Innovative ideas presented in the first two chapters are mainly related to 
the ways in which difficult subjects can be taught to beginner students.
The first really remarkable scientific concept which appears
in this book is, in our opinion, presented in Chapter III. The author begins 
with a short and rather classical account of the general principles of classical mechanics.
Then he proves that on the evolution space of a dynamical system, there exists
a remarkable presymplectic form, which he calls the \emph{Lagrange form}.
Its kernel determines the vector field whose integral curves are the 
\emph{motions} of the system. This form projects onto the set of all possible
motions, called the \emph{space of motions}, and its
projection is a \emph{symplectic form}, \emph{i.e.}, a closed, non degenerate $2$-form.
The author chooses this property, called the \emph{Maxwell principle}, as 
the founding principle of mechanics. For us, it is a very important innovation:
the traditionally used concepts, such as the \emph{configuration space}, the
\emph{evolution space}, the \emph{phase space}, go into the background and 
leave the first place to the \emph{space of motions} and to the \emph{Lagrange form} 
of which it is endowed.  Thanks to this innovation, the author will be able to describe,
in the framework of classical mechanics, with the dynamics of a particle 
with spin, though there is no Lagrangian for such a system, and with assemblies of 
an indeterminate number of such particles. Similarly, he will be able to 
describe systems made by massless relativistic particles, though there is no
evolution space for such systems.
\par\smallskip

The space of motions of a dynamical system is not very often considered
by modern authors, though it appeared as soon as 1808 in the works of Lagrange. 
This very natural concept has a nice mathematical property: the space of motions is 
always endowed with a smooth manifold structure. 
The \emph{phase portrait} of a dynamical system, a frequently used closely related concept,  
very often has a much more complicated structure. One may wonder why
the concept of space of motions is not used more by modern authors. 
Maybe it is because for some dynamical systems, the space of motions 
is a \emph{non-Hausdorff} manifold. Another possible explanation is that
some scientists are interested in the thorough description of particular motions of a system,
rather than by the study of the set of all possible motions. 
By showing that many important results can be deduced from the symmetries
of the space of motions of a system, the author proves that this reluctance is ungrounded.
For example, in Chapter III, section 12, he proves that the principle of 
equality of action and reaction appears as a consequence of Galilean
relativity and Maxwell's principle.
\par\smallskip

In classical (non-relativistic) mechanics, the cohomological interpretation 
of the \emph{mass} of an isolated dynamical system  is, in  our opinion,
another innovative idea worth mentioning. A cohomology class appears indeed
in the mathematical expression of the action of the Galilei group on the space of
motions of such a system. It is an element of a one-dimensional vector space,
and the author has proven that it can be interpreted as the mass of the system.
In classical mechanics, it is legitimate to consider isolated systems with 
a positive, null or negative mass. It will be done by the author in his study
of the behaviour of elementary systems with respect to time or space reversals. 
\par\smallskip

Another originality of this book is that it presents dynamical systems 
in the framework of relativistic physics as well as in the framework of 
classical, non-relativistic mechanics. This is made possible thanks to 
the use of the concept of state of motions. In classical, non-relativistic mechanics,
the Kepler problem and the $n$-bodies problem are very good mathematical 
models for systems of gravitationally interacting material points.
In relativistic physics, there is no such mathematical model for a system of
interacting electrically charged  particles.  These particles interact 
by means of the electomagnetic field created by their motion. Once created, this 
electromagnetic field evolves according to its own laws. Up to now, no mathematical model
is available to  describe the motions of the particles with the integral curves 
of a vector field defined on some hypothetic evolution space, only depending on the
positions and motions of the particles. The author manages to get around this difficulty
by using the space of motions of the system and Maxwell's principle.
\par\smallskip

In classical mechanics, the symmetry group of the space of motions of an 
isolated dynamical system is the Galilei group. In relativistic physics, 
this group is the Poincaré group. This change of symmetry group has 
important consequences clearly described by the author.
The status of \emph{mass} is in relativistic physics very different from 
its status in classical mechanics, due to the fact that the symplectic cohomology 
of the Poincaré group is trivial. In relativistic physics, there is no more
a \emph{barycentric decomposition} of the state of motions of an isolated system,
as in classical mechanics, because the Poincaré group has no privileged
normal Abelian subgroup, as the Galilean group has.
\par\smallskip

The works of the German mathematician Emmy Noether \cite{Noether} told us that 
\emph{first integrals} of a Lagrangian or Hamiltonian system very often
are linked with symmetries of the equations of motion. 
Most textbooks in classical mechanics published before
the author's book only indicate how a real valued first integral is 
determined by each one-parameter symmetry group of these equations.
We believe that Jean-Marie Souriau is, with Stephen Smale \cite{Smale}, 
among the first scientists who considered the geometric structure of the
full set of these first integrals: they are the components of the \emph{moment map}
of the symmetry group's action, defined on the evolution space and valued in the dual
vector space of the Lie algebra of this group. 
\par\smallskip

An \emph{elementary system} is a relativistic dynamical system whose 
space of motions is an homogeneous space of the restricted Poincaré group.
This purely mathematical definition, due to the author, seems to us very important 
because the classification of elementary systems reveals many properties of 
physicits' \emph{elementary particles}, especially their \emph{mass}, their
\emph{spin} and, for massless particles, their \emph{helicity} (whose 
possible values, $1$ or $-1$, correspond to the two different circular 
polarizations of photons). In the non-relativistic approximation, 
elementary systems are models of elementary particles in the 
framework of classical mechanics. By this means, Geometric Optics appears, when
the particle's spin is negligible, as a part of classical mechanics. 
The author's works so appear in continuity with those of Hamilton \cite{Hamilton1},
who introduced the \emph{characteristic function} first in Optics before 
using it in mechanics. Interested readers will find in the long Introduction of the
book \cite{GuilleminSternberg84b} a very nice discussion of the symplectic aspects
of geometric Optics and Electromagnetism. 
\par\smallskip

Chapter IV also contains several worth mentioning innovations. For the author,
a \emph{statistical state} of a dynamical system is a probability measure 
on the space of motions of the system. This definition is more natural than
that generally used, for example by George~Mackey  \cite{Mackey1963} 
who, instead of the state of motions, uses the phase space at a given time.
However, the most important innovation contained in this chapter seems to us
the generalization of the notion of a Gibbs state, in which the energy is 
replaced by the moment map or the action of a symmetry group. This very 
natural generalization (the energy being the moment 
map of the action of the one-dimensional group of time translations) involves remarkable new features 
when the symmetry group is not Abelian. This group acts on the dual space 
of its Lie algebra by an affine action involving a symplectic cocycle,
with the coadjoint action as linear part, for which the moment map is equivariant.   
By using generalized Gibbs states, the author develops a kind of
\emph{Lie groups thermodynamics}, which seems to us very interesting, as well
from the mathematician's viewpoint as for possible applications in physics.
In this theory, the generalized temperature and the generalized quantity of heat
are, respectively, elements of the Lie algebra of the symmetry group and of its
dual vector space. The author proves that an open convex subset of the Lie algebra
is endowed with a Riemannian metric which plays an important part in Statistics and
Information theory.
\par\smallskip

The seeds of Geometric Quantization, presented in Chapter V, can be found 
in George Mackey's small book \cite{Mackey1963}, published in 1963. 
Jean-Marie Souriau \cite{SSD} and Bertram Kostant \cite{Kostant70}, 
its  main creators, independently proposed two slightly
different, but equivalent versions of this theory. On the basis made by 
a symplectic manifold, Kostant defines a bundle whose fibres are complex lines,
endowed with a connection whose curvature is the symplectic form of its basis.
Similarly, Souriau defines a bundle whose fibres are circles, endowed with 
a contact form $\omega$, whose Reeb vector field is tangent to the fibres.
Moreover, the exterior derivative of $\omega$ projects onto the symplectic 
form of the basis. Souriau's circle bundle is the \emph{principal bundle}
which is associated to Kostant's complex line bundle, and the connection form 
of Kostant's line bundle is the contact form $\omega$ of Souriau.     
While Souriau quantizes the \emph{manifold of motions} of a classical
dynamical system, Kostant quantizes its \emph{phase space} at a given time.
However, this difference in their choice of the quantized symplectic manifold 
may not be as important as it seems, because Souriau often uses the local 
symplectomorphism which exists between the space of motions and the phase space, 
especially when he derives the  quantized wave equations.
\par\smallskip

Many remarkable results are presented in Chapter V. For instance, the author
proves that a quantizable system made by indistinguishable particles has exactly
two non-equivalent quantizations, both experimentally observed: they are the 
Fermi-Dirac and Bose-Einstein quantizations well known by physicists. He proves too
that wave functions of quantum mechanics can be described with functions defined on the
quantum manifold, and that the usual spatio-temporal description is obtained 
by means of a Fourier transform. The quantization of the space of motions
of an assembly of indistinguishable particles leads to the \emph{second 
quantization formalism}: quantum vacuum, creation and annihilation operators.
The author also shows that the geometric obstructions encountered when the action
of a symmetry group on the classical space of motions of a system is 
lifted to the quantum manifold, explains some facts well known by physicists: 
the lift of the Abelian group of translations is the non-Abelian 
\emph{Weyl group}; elements of the Galilean group can be separately lifted, 
but this full group does not act on the quantum manifold (fact due to a 
cohomological obstruction discovered by Valentine Bargmann \cite{bargmann});
the Poincaré group can be lifted for a particle of integer spin, but not
for a particle of half-integer spin, in which case it is a
two sheets covering space of the Poincaré group which can be lifted, which 
leads to the use of Dirac's spinors.
\par\smallskip

Finally, we must quote the many conjectural ideas presented by the author
in the last section of his book, for example about the behaviour of elementary 
systems with respect to space or time reversals, and about a generalization of
the notion of quantum state.

\subsection*{Language and notations of the author} 
The author's language, although slightly different from the one generally used 
in differential geometry, is perfectly logical and understandable. One of 
the very few terms which could cause misunderstandings is that of 
\emph{embedding} (in French, \emph{plongement}). For the author, an embedding is
a smooth map defined on a smooth manifold $V$, or on an open subset of $V$, 
with values in another smooth manifold $V'$, which is injective (the author uses 
the term \emph{regular} for \emph{injective}) and whose rank is everywhere equal
to the dimension of $V$. Most geometers call such a map 
\emph{an injective immersion}, and use the term \emph{embedding} for injective
immersions which, in addition, are homeomorphisms of their domain of definition 
onto their image, endowed with the topology induced by that of $V'$. It seems 
to us that the choice made by the author is very reasonable because injective 
immersions are much more frequently encountered than embeddings. For example,
orbits of a Lie group action, as well as leaves of a foliation, always are 
immersed in the manifold in which they are contained, and much more rarely embedded.
This indication only concerns the original version of the book in French. 
In its English translation \cite{SSDEng}, the translators use the usual term
\emph{injective immersion} and, in a footnote, indicate that for the author a 
\emph{submanifold} is an \emph{immersed submanifold}.
\par\smallskip 

Another particularity of the author's language which seems to us worth mentioning,
is about the notion of an \emph{Euclidean vector space} of finite dimension. 
For the author, it is a real or complex finite-dimensional vector space $E$
endowed with a symmetric, bilinear form $g$ satisfying the 
following conditions.
\par\smallskip 

\begin{itemize}

\item{} If $E$ is a real vector space, $g$ is assumed to be non degenerate, which 
means that for any $x\in E$, $x\neq 0$, there exits $y\in E$ such that $g(x,y)\neq 0$.
The author does not impose the condition $g$ positive.
He still uses the term \emph{Euclidean vector space}, not the term 
\emph{pseudo-Euclidean vector space} when $g$ is definite without being positive, 
as done by most mathematicians.

\item{} If $E$ is a complex $n$-dimensional vector space, $g$ is a symmetric, 
bilinear form \emph{for the structure of $2n$-dimensional real vector space of}
$E$ underlying its structure of $n$-dimensional complex vector space. The form $g$
is assumed to be non degenerate and to satisfy, for all pair $(x,y)$ of elements 
in $E$, $g(ix,iy)=g(x,y)$, where $i=\sqrt{-1}$. 
\end{itemize}
This convention seems to us very logical and useful, because it unifies two different
concepts often taught separately: the concept of real Euclidean finite-dimensional
vector space and that of complex finite-dimensional Hermitian vector space.  
\par\smallskip

Some notations of the author seem to us disconcerting. For example, he calls
\emph{variable} any symbol, such as the letter $y$, and associates it to a map. 
He denotes by $[\hbox{value}\ y](a)$ the value of the map associated to $y$ at the point
$a$ if its domain of definition. He uses different symbols for the variable, 
the associated map, the point of its domain of definition and its value at that point. 
So it is sometimes difficult to understand some expressions he writes, because the reader
must simultaneously have in mind the meaning of many symbols.
\par\smallskip

Given a vector field $f$ defined on a smooth manifold $V$, 
the author calls \emph{derivation} associated to $f$ the operation which, 
to any smooth map $A$ defined on $V$ and valued in a smooth manifold $V'$, 
associates the map whose value, at each $x\in V$, is the image of the vector 
$f(x)$ by the linear map tangent to $A$ at $x$. He denotes derivations 
by symbols such as $d$ or $\delta$, other than the symbol which denotes 
the vector field (here the symbol which denotes the vector field is $f$).
When the author considers several vector fields, he uses symbols
such as $\d$, $\d_1$, $\delta$, $\delta_1$, etc to denote the corresponding 
derivations. Of course, these conventions are perfectly logical and coherent.
However, their use may make the reading of some formulae rather difficult.
Moreover, since the symbol $\d$ is used for derivations associated to vector fields, 
the author cannot use it for the exterior derivation of differential forms, as done
by most mathematicians. For the exterior derivation of differential forms, the 
author uses the symbol $\nabla$, which may disconcert some readers. 
\par\smallskip

Let us finally indicate a few pecularities of the notations used by the author in
exterior differential calculus. He does not use the symbol $\wedge$ 
for the exterior product of differential forms, which makes some formulae heavy.
He denotes by $\eta(X)$ the interior product of a differential form $\eta$
by a vector field $X$, while many other mathematicians denote it by
$\mathi(X)\eta$,  $\mathi_X\eta$ or $\mathi_X(\eta)$. It is probably for this 
reason that he denotes by $\eta(X_1)(X_2)\cdots(X_p)$ 
the evaluation of the $p$-form $\eta$ on the vector fields $X_1,X_2,\ldots,X_p$, 
instead of writing more simply $\eta(X_1,\ldots,X_p)$.  At first glance, 
this convention may seem simple and convenient. However, it makes some 
formulae heavy, such as the Cartan formula 
${\mathcal L}(X)=\mathi(X)\d+\d \mathi(X)$,
which expresses the Lie derivative with respect to a vector field $X$ 
as the anticommutator of the interior product by $X$ with the exterior derivative.

\section{Novel research ideas in Jean-Marie Souriau's footsteps}

Every research book is a survey at a given time of the state of knowledge 
obviously limited on the issue. Just as the author considered that the book by Lagrange 
\textit{Mécanique analytique} was unfinished, we believe that \textit{Structure des 
systèmes dynamique} is a very in-depth work opening new research paths. 

\section{Conclusion}

We cannot fail to be impressed when reading this book by the extent and the thoroughness 
of the author's knowledge, as well in mathematics as in mechanics or in physics, 
and by the originality and the depth of his thoughts. Jean-Marie Souriau is the author of two other
very remarkable books, \emph{Géométrie et Relativité} \cite{GR}
and \emph{Calcul linéaire} \cite{CL}, which are very rich and original too
and deserve to be read, and read again. We believe that among the paths for research 
he indicates in these books, many still are not yet fully explored.

\end{document}